\documentclass[11pt,a4paper]{article}
\usepackage{amsfonts}
\usepackage{}

\usepackage{graphicx}
\usepackage{amsfonts,amsmath, amssymb}

\newtheorem{theo}{\textbf{\ \ \quad Theorem}}[section]

\newtheorem{lem}{\textbf{\ \ \quad Lemma}}[section]
\newtheorem{remark}{\textbf{\ \ \quad Remark}}[section]

\newtheorem{prop}{\textbf{\ \ \quad Proposition}}[section]

\newcommand{\lbl}[1]{\label{#1}}

\newcommand{\be}{\begin{equation}}
\newcommand{\ee}{\end{equation}}
\newcommand\bes{\begin{eqnarray}}
\newcommand\ees{\end{eqnarray}}
\newcommand{\bess}{\begin{eqnarray*}}
\newcommand{\eess}{\end{eqnarray*}}

\newcommand{\nm}{\nonumber}

\newcommand{\ds}{\displaystyle}

{\setlength\arraycolsep{2pt}}

\title{Global existence and non-existence of stochastic parabolic equations }

\author{Guangying Lv$^{a,b}$,  Jinlong Wei$^c$\\
\\
\ \\
   {\small \it $^a$ Institute of Applied Mathematics, Henan University, Kaifeng, Henan 475001, China}\\
  {\small \it $^b$ Center for Applied Mathematics, Tianjin University, Tianjin 300072, China }\\
   {\small \tt gylvmaths@henu.edu.cn}\\
   {\small \it $^c$ School of Statistics and Mathematics, Zhongnan University of}\\
   {\small \it Economics and Law, Wuhan 430073, China}\\
    {\small \tt  weijinlong.hust@gmail.com }
}

\begin{document}
\maketitle

\medskip

\begin{abstract}
This paper is concerned with the blowup phenomenon of stochastic parabolic equations both on
bounded domain and in the whole space.
We introduce a new method to study the blowup phenomenon on bounded domain.
Comparing with the existing results, we delete the assumption that the solutions to
stochastic heat equations are non-negative. Then the blowup phenomenon in the whole space is
obtained by using the properties of heat kernel. We obtain that
the solutions will blow up in finite time for nontrivial initial data.



\end{abstract}

\baselineskip=15pt

\section{Introduction}
\setcounter{equation}{0}

For deterministic partial differential equations,
finite time blowup phenomenon has been studied by many authors, see the book \cite{Hubook2018}.
There are two cases to study this problem. One is bounded domain and the other is whole space.
On the bounded domain, the $L^p$-norm of solutions ($p>1$) will blow up in finite time.
The methods used for bounded domain include:  Kaplan's first eigenvalue method, concavity method and
comparison method, see Chapter 5 of \cite{Hubook2018}. The main result is the following: under the assumptions that
the initial data is suitable large and that the nonlinear term $f(u)$ satisfies $f(u)\geq u^{1+\alpha}$
with $\alpha>0$, the solution of $u_t-\Delta u=f(u)$ with Dirichlet boundary condition will blow up in finite time.

For the whole space, the following "Fujita Phenomenon" has been attraction
in the literature. Consider the following Cauchy problem
   \bes\left\{\begin{array}{lll}
u_t=\Delta u+u^p,\ \ \ &x\in\mathbb{R}^d, \ \ t>0,\ \ p>0,\\
u(0,x)=u_0(x), \ \ \ &x\in\mathbb{R}^d.
   \end{array}\right.\lbl{1.1}\ees
It has been proved that:
 \begin{quote}
(i) if $0<p<1$, then every nonnegative solution is global, but not necessarily
unique;

(ii) if $1<p\leq1+\frac{2}{d}$, then any nontrivial, nonnegative solution
blows up in finite time;

(iii) if $p>1+\frac{2}{d}$, then $u_0\in\mathcal{U}$ implies that
$u(t,x,u_0)$ exists globally;

(iv) if $p>1+\frac{2}{d}$, then $u_0\in\mathcal{U_1}$ implies that
$u(t,x,u_0)$ blows up in finite time,
 \end{quote}
where $\mathcal{U}$ and $\mathcal{U_1}$ are defined as follows
   \bess
\mathcal{U}&=&\left\{v(x)|v(x)\in BC(\mathbb{R}^d,\mathbb{R}_+),
v(x)\leq \delta e^{-k|x|^2},\ k>0,\delta=\delta(k)>0\right\},\\
\mathcal{U_1}&=&\left\{v(x)|v(x)\in BC(\mathbb{R}^d,\mathbb{R}_+),
v(x)\geq c e^{-k|x|^2},\ k>0,c\gg1\right\}.
  \eess
Here $BC=\{$ bounded and uniformly continuous functions $\}$, see
Fujita \cite{F1966,F1970} and Hayakawa \cite{kH1973}.

It is easy to see that for the whole space, there are four types of behaviours for problem (\ref{1.1}), namely,
(1) global existence unconditionally but uniqueness fails in certain solutions, (2) global existence with restricted
initial data, (3) blowing up unconditionally, and (4) blowing up with restricted initial data. The
occurrence of these behaviors depends on the combination effect of the nonlinearity represented by the parameter $p$, the size of the initial datum $u_0(x)$, represented by the choice of $\mathcal{U}$
or $\mathcal{U_1}$, and the dimension of the space.

Now, we recall some known results of stochastic partial differential equations (SPDEs).
In this paper, we only focus on the stochastic parabolic equations.
It is known that the existence and uniqueness of
global solutions to SPDEs can  be   established under appropriate conditions (\cite{Cb2007,DW2014,LR2010,L2013,T2009}).
For the finite time blowup phenomenon of stochastic parabolic equations, we first consider
the case on bounded domain. Consider the following equation
   \bes\left\{\begin{array}{lll}
   du=(\Delta u+f(u))dt+\sigma(u)dW_t, \ \ \qquad t>0,&x\in  D,\\[1.5mm]
   u(x,0)=u_0(x)\geq0, \ \ \ &x\in D,\\[1.5mm]
   u(x,t)=0, \qquad \qquad \qquad \qquad \qquad \qquad t>0,  &x\in\partial D,
    \end{array}\right.\lbl{1.2}\ees
Da Prato-Zabczyk \cite{PZ1992} considered  the existence of global solutions of (\ref{1.2})  with
additive noise ($\sigma$ is constant). Manthey-Zausinger \cite{MZ1999} considered (\ref{1.2}),
where $\sigma$ satisfied the global Lipschitz condition. Dozzi and L\'{o}pez-Mimbela
\cite{DL2010} studied   equation (\ref{1.2}) with
$\sigma(u)=u$ and  proved that if $f(u)\geq u^{1+\alpha}$ ($\alpha>0$) and initial data is large enough, the solution will blow up in finite time,
and that if $f(u)\leq u^{1+\beta}$ ($\beta$ is a certain positive constant) and the initial data is
small enough, the solution will exist globally, also see \cite{NX2012}. A natural question arises: If
$\sigma$ does not satisfy the global Lipschitz condition, what can we say about the solution?
Will it blow up in finite time or exist globally? Chow \cite{C2009,C2011} answered part of this question.
Lv-Duan \cite{LD2015} described the competition between the nonlinear term and noise term for equation (\ref{1.2}).
Bao-Yuan \cite{BY2014} and Li et al.\cite{LPJ2016} obtained the existence of local solutions of (\ref{1.2}) with jump process
and L$\acute{e}$vy process, respectively. For blowup phenomenon of stochastic functional
parabolic equations, see \cite{CL2012,LWW2016} for details.
 In a somewhat different case, Mueller \cite{M1991} and, later,
Mueller-Sowers \cite{MuS1993} investigated the problem of a noise-induced explosion for a
special case of equation (\ref{1.2}), where $f(u)\equiv0,\,\sigma(u)=u^\gamma$ with
$\gamma>0$ and $W(x,t)$ is a space-time white noise. It was shown that the solution will
explode in finite time with positive probability for some $\gamma>3/2$.

We remark that the method used to prove the finite time blowup on bounded domain is
the stochastic Kaplan's first eigenvalue method. In order to make sure the inner product
$(u,\phi)$ is positive, the authors firstly proved the solutions of (\ref{1.2}) keep positive under some
assumptions, see \cite{BY2014,C2009,C2011,LPJ2016,LD2015}.
We find that under some special case the positivity of solution can be deleted.
What's more, in present paper, we will give
a new method (stochastic concavity method) to prove the solutions blow up in finite time.
The advantage of this method is that we need not the positivity of solution.

For the whole space, Foondun et al. \cite{FLN2018} considered the finite time blowup phenomenon
for the Cauchy problem of stochastic parabolic equations.
Comparing with the deterministic parabolic equations, they only obtained
the result similar to type (4).
In this paper, we establish the similar results to types (1)
and (3). The method used here is comparison principle and the properties of heat kernel.
We obtain some different phenomenon with or without noise. Moreover,
many types of noise are considered.

Comparing with the results of deterministic partial differential equations,
there are a lot of work to do and we will study this issue in our further paper.

This paper is arranged as follows.   In Sections 2 and 3,
we will consider the global existence and non-existence of stochastic parabolic equations
on bounded domain and in the whole space, respectively.
This paper ends with a short discussion in Section 4.

Throughout this paper, we write $C$ as a general positive constant and $C_i$, $i=1,2,\cdots$ as
a concrete positive constant.

\section{Bounded domain}
\setcounter{equation}{0}
In this section, we first recall some known results on bounded domain, and then
give some non-trivial generalizations. Consider  the following SPDE
      \bes\left\{\begin{array}{lll}
   du=(\Delta u+f(u,x,t))dt+\sigma(u,\nabla u,x,t)dW_t, \ \ \qquad &t>0,\ x\in  D,\\[1.5mm]
   u(x,0)=u_0(x), \ \ \ &\qquad\ \ x\in D,\\[1.5mm]
   u(x,t)=0, \ \ \ \ &t>0,\ x\in\partial D,
    \end{array}\right.\lbl{2.1}\ees
where $\sigma$ is a given function, and $W(x,t)$ is a Wiener random field defined in a complete probability
space $(\Omega,\mathcal {F},\mathbb{P})$ with a filtration $\mathcal {F}_t$. The Wiener random field has mean
$\mathbb{E}W(x,t)=0$ and its covariance function $q(x,y)$ is defined by
    \bess
\mathbb{E}W(x,t)W(y,s)=(t\wedge s)q(x,y), \ \ \ x,y\in\mathbb{R}^n,
   \eess
where $(t\wedge s)=\min\{t,s\}$ for $0\leq t,s\leq T$. The existence
of strong solutions of (\ref{2.1}) has been studied by many authors  \cite{Cb2007,PZ1992}.
To consider positive solutions, they  start with  the unique
solution $u\in C(\bar D\times[0,T])\cap L^2((0,T);H^2)$ for equation (\ref{2.1}).
Chow \cite{C2009,C2011} considered the finite time blowup problem of (\ref{2.1}).
They used the positivity of solution to prove the finite time blowup.
Under the following conditions
   \begin{quote}
 (P1) There exists a constant $\delta\geq0$ such that
 \bess
 \frac{1}{2} q(x,x)\sigma^2(r,\xi,x,t)-\sum_{i,j=1}^na_{ij}\xi_i\xi_j\leq\delta r^2
    \eess
 for all $r\in\mathbb{R},x\in\bar D,\xi
\in\mathbb{R}^n$ and $t\in[0,T]$;\\
 (P2) The function $f(r,x,t)$ is continuous on $\mathbb{R}\times\bar D\times[0,T]$ and such that
 $f(r,x,t)\geq0$ for $r\leq0$ and $x\in\bar D$, $t\in[0,T]$;  and \\
 (P3) The initial datum $u_0(x)$ on $\bar D$ is positive and continuous,
 \end{quote}

\begin{prop}\lbl{p1.1}{\rm\cite[Theorem 3.3]{C2009}} Suppose that the conditions
{\rm (P1),(P2)} and {\rm (P3)} hold true. Then the solution of the initial-boundary problem
for the parabolic It\^{o}'s equation {\rm(\ref{2.1})} remains positive, i.e.,
$u(x,t)\geq 0$, a.s. for almost every $x\in D$ and for all $ t\in[0,T]$.
\end{prop}

Let $\phi$ be the eigenfunction with respect to the first eigenvalue $\lambda_1$ on the bounded
domain, i.e.,
  \bess\left\{\begin{array}{llll}
-\Delta \phi=\lambda_1 \phi, \ \ \ \ \ \ \  \ {\rm in} \  D,\\
\phi=0, \ \ \qquad\ \qquad  {\rm on}\ \partial D.
   \end{array}\right.\eess
And we   normalize it in such a way that
   \bess
\phi(x)\geq0,\ \ \ \ \int_ D \phi(x)dx=1.
   \eess
In paper \cite{C2011}, Chow assumed that the following conditions hold

   \begin{quote}
   (N1) There exist a continuous function $F(r)$ and a constant $r_1>0$ such that
   $F$ is positive, convex and strictly increasing for $r\geq r_1$ and satisfies
      \bess
f(r,x,t)\geq F(r)
   \eess
for $r\geq r_1$, $x\in\bar D$, $t\in[0,\infty)$;\\
(N2) There exists a constant $M_1>r_1$ such that $F(r)>\lambda_1r$ for $r\geq M_1$;\\
(N3) The positive initial datum satisfies the condition
   \bess
(\phi,u_0)=\int_ D u_0(x)\phi(x)dx>M_1;
   \eess
(N4) The following condition holds
   \bess
\int_{M_1}^\infty\frac{dr}{F(r)-\lambda_1r}<\infty.
   \eess
    \end{quote}
Alternatively, he imposes the following conditions $S$ on the noise term:
   \begin{quote}
(S1) The correlation function $q(x,y)$ is continuous and positive for $x,y\in\bar D$
such that
   \bess
\int_ D\int_ D q(x,y)v(x)v(y)dxdy\geq q_1\int_ D v^2(x)dx
   \eess
for any positive $v\in H$ and for some $q_1>0$;

(S2) There exist a positive constant $r_2$, continuous functions $\sigma_0(r)$ and $G(r)$ such that
they are both positive, convex and strictly increasing for $r\geq r_2$ and satisfy
  \bess
\sigma(r,x,t)\geq \sigma_0(r)\ \ \ \ {\rm and} \ \ \ \ \sigma_0^2(r)\geq2G(r^2)
   \eess
for $x\in\bar D$, $t\in[0,\infty)$;

(S3) There exists a constant $M_2>r_2$ such that $q_1G(r)>\lambda_1r$ for $r\geq M_2$;

(S4) The positive initial datum satisfies the condition
   \bess
(\phi,u_0)=\int_ D u_0(x)\phi(x)dx>M_2;
   \eess

(S5) The following integral is convergent so that
   \bess
\int_{M_2}^\infty\frac{dr}{q_1G(r)-\lambda_1r}<\infty.
   \eess
    \end{quote}

\begin{prop}\lbl{p1.2} {\rm\cite[Theorem 3.1]{C2011}}
Suppose the initial-boundary value problem {\rm(\ref{2.1})} has a unique local solution
and the conditions {\rm(P1)-(P3)} are satisfied, where $\sigma$ does not depend on
$\nabla u$. In addition, we assume that either
the conditions {\rm(N1)-(N4)} or the alternative conditions {\rm(S1)-(S5)} given above
hold true. Then, for a real number $p>0$, there exists a constant $T_p>0$ such that
   \bess
\lim\limits_{t\rightarrow T_p-}\mathbb{E}\|u \|_p
=\lim\limits_{t\rightarrow T_p-}\mathbb{E}\left(\int_ D|u(x,t)|^pdx\right)^\frac{1}{p}=\infty,
   \eess
where $p\geq1$ under conditions $N$, while $p\geq2$ under conditions $S$.
\end{prop}

The positivity of solutions is needed for the case that the nonlinear term induces the
finite time blowup.
But for a special case, we can prove the positivity of solutions can be deleted. Now, we consider
the following SPDEs
      \bes\left\{\begin{array}{lll}
   du=\Delta udt+\sigma(u,x,t)dW(x,t), \ \ \qquad &t>0,\ x\in  D,\\[1.5mm]
   u(x,0)=u_0(x), \ \ \ &\qquad\ \ x\in D,\\[1.5mm]
   u(x,t)=0, \ \ \ \ &t>0,\ x\in\partial D,
    \end{array}\right.\lbl{2.3}\ees
where $W(x,t)$ is time-space white noise and $D\subset\mathbb{R}$ is an interval in $\mathbb{R}$.
    \begin{theo}\lbl{t2.1} Assume that the initial-boundary problem (\ref{2.3}) has
a unique local solution.
Assume further that $C_1|u|^\gamma\leq|\sigma(u,x,t)|\leq C_2|u|^\gamma_1$ with $C_1>0$ and $\gamma_1\geq\gamma>1$, $u_0\geq0$ and
   \bess
\left(\int_ D u_0(x)\phi(x)dx\right)^{2(\gamma-1)}\geq\frac{\lambda_1}{q_1C_1^2}.
   \eess
Then there exist constants $T^*>0$ and $p\geq2\gamma_1$ such that
     \bess
\lim\limits_{t\rightarrow T^*-}\mathbb{E}\|u_t\|^p_{L^p}
=\lim\limits_{t\rightarrow T^*-}\mathbb{E}\int_D|u(x,t)|^pdx=\infty.
     \eess
  \end{theo}

{\bf Proof.}
 We will prove the theorem by contradiction.
Suppose finite time blowup is false. Then there exist a global positive solution $u$ and $p\geq2\gamma_1$ such
that for any $T>0$
   \bess
\sup_{0\leq t\leq T}\mathbb{E}\|u(\cdot,t)\|^p_{L^p}<\infty,
   \eess
which implies that
   \bess
\sup_{0\leq t\leq T}\mathbb{E}\Big|\int_ D u(x,t)\phi(x)dx\Big|^2\leq\|\phi\|^2_{L^q(D)}\sup_{0\leq t\leq T}\mathbb{E}\|u(\cdot,t)\|^p_{L^p}<\infty,
    \eess
where $1/p+1/q=1$, $\phi$ is defined as below Proposition \ref{p1.1} and satisfies $\int_ D\phi(x)dx=1$.
Define
   \bess
\hat u(t):=\int_ D u(x,t)\phi(x)dx.
   \eess
By applying It\^{o}'s formula to $\hat u^2(t)$, we get
   \bes
\hat u^2(t)&=&(u_0,\phi)^2-2\lambda_1\int_0^t\hat u^2(s)ds+2\int_0^t\int_ D \hat u(s)\sigma(u,x,t)\phi(x)d  W(x,s)dx\nm\\
&&+\int_0^t\int_ D\sigma^2(u,x,s)\phi^2(x)dxds
   \lbl{2.4}\ees
We note that the stochastic term is usually a local martingale. Thus we need use the technique of stopping time.
Let
   \bess
\tau_n=\inf\{t\geq0:\ \int_0^t\int_ D\sigma^2(u,x,s)\phi^2(x)dxds\geq n\}.
  \eess
Let $\eta(t\wedge\tau_n)=\mathbb{E}\hat u^2(t\wedge\tau_n)$. By taking an expectation over (\ref{2.4}), we obtain
  \bess
\eta(t\wedge\tau_n)&=&(u_0,\phi)^2-2\lambda_1\int_0^{t\wedge\tau_n}\eta(s)ds
+\int_0^{t\wedge\tau_n}\mathbb{E}\int_ D\sigma^2(u,x,s)\phi^2(x)dxds.
    \eess
Noting that
   \bess
\eta(t\wedge\tau_n)\leq(u_0,\phi)^2
+\int_0^{t}\mathbb{E}\int_ D\sigma^2(u,x,s)\phi^2(x)dxds,
    \eess
and letting $n\to\infty$, we have
  \bess
\eta(t)=(u_0,\phi)^2-2\lambda_1\int_0^{t}\eta(s)ds
+\int_0^{t}\mathbb{E}\int_ D\sigma^2(u,x,s)\phi^2(x)dxds.
    \eess
Using the assumptions $\inf_{x,y\in D}q(x,y)\geq q_1>0$ and $\sigma^2(u,x,s)\geq C_1|u|^{2\gamma}$ with $\gamma>1$ and
Jensen's inequality, we have
   \bess
\eta(t)&\geq&\eta(0)-2\lambda_1\int_0^t\eta(s)ds
+2q_1C_1^2\int_0^t\eta^\gamma(s)ds,
   \eess
or, in the differential form,
   \bess\left\{\begin{array}{llll}
\ds\frac{d\eta(t)}{dt}=-2\lambda_1\eta(t)+2q_1C_1^2 \eta^\gamma(t)\\[1.5mm]
\eta(0)=\eta_0.
    \end{array}\right. \eess
Noting that
   \bess
\eta(0)=\left(\int_ D u_0(x)\phi(x)dx\right)^{2}\geq\left(\frac{\lambda_1}{q_1C_1^2}\right)^{\frac{1}{(\gamma-1)}},
   \eess
we have $\eta'(0)\geq0$. This implies that $\eta(t)>0$. An integration of the differential equation gives that
   \bess
T\leq\int_{\eta_0}^{\eta(T)}\frac{dr}{C^2_1q_1r^{\gamma}-\lambda_1r}
\leq\int_{\eta_0}^\infty\frac{dr}{C^2_1q_1r^{\gamma}-\lambda_1r}<\infty,
    \eess
which implies $\eta(t)$ must blow up at a time $T^*\leq\int_{\eta_0}^\infty\frac{dr}{C^2_1q_1r^{\gamma}-\lambda_1r}$.
Hence this is a contradiction. This completes the proof. $\Box$

The advantage of Theorem \ref{t2.1} is that the positivity of the solution is not needed.
And in above Theorem, we assume that the initial-boundary problem (\ref{2.3}) has
a unique local solution. In fact, if $\sigma$ satisfies the local Lipschitz condition,
one can follow the method of \cite{T2009} to obtain the existence and uniqueness of local solution, also see \cite{eP1979}.
In \cite{eP1979,T2009}, the authors established the existence and uniqueness of energy solution, where
the solutions belong to $H_0^1(D)$ for any fixed time almost surely. Noting that
$H^{\frac{1}{2}+}(D)\hookrightarrow L^\infty(D)$ for $D\subset\mathbb{R}$, our assumptions are valid.

If we only consider the case
$\sigma$ does not depend on $\xi$, that is,
$\sigma:=\sigma(u,x,t)$.  Then it follows the assumption (P1) that
$\sigma(0,x,t)=0$, which implies that for additive noise, the solutions maybe not
keep positive. Hence the first eigenvalue method will fail. Next, we introduce another method.
For simplicity, we consider the following SPDEs
       \bes\left\{\begin{array}{lll}
   du=[\Delta u+|u|^{p-1}u]dt+\sigma(x,t)dB_t, \ \ \qquad &t>0,\ x\in  D,\\[1.5mm]
   u(x,0)=u_0(x), \ \ \ &\qquad\ \ x\in D,\\[1.5mm]
   u(x,t)=0, \ \ \ \ &t>0,\ x\in\partial D,
    \end{array}\right.\lbl{2.5}\ees
where $B_t$ is an one-dimensional Brownian motion. If the initial data belongs to $H^1(D)$,
Debussche et al. \cite{DMH2015} proved the solution of (\ref{2.5}) belongs to $H^3_0(D)$ during the
lifespan.
   \begin{theo}\lbl{t2.2} Suppose that $p>1$ and $u_0$ satisfies
      \bess
-\frac{1}{2}\int_D|\nabla u_0(x)|^2dx+\frac{1}{p+1}\int_D|u_0(x)|^{p+1}dx-\frac{1}{2}\int_0^\infty\mathbb{E}\int_D|\nabla\sigma(x,t)|^2dxdt>0,
   \eess
then the solution of (\ref{2.5}) must blow up in finite time in sense of mean square.
  \end{theo}

{\bf Proof.} We will
prove the theorem by contradiction.
First we suppose there exist a global solution $u$ such that
    \bess
\sup_{t\in[0,T]}\mathbb{E}\int_Du^2dx<\infty
   \eess
for any $T>0$. Similar to the proof of Theorem \ref{t2.1}, by using It\^{o} formula, we have
  \bess
\mathbb{E}\int_Du^2-\int_Du^2_0=-2\mathbb{E}\int_0^t\int_D|\nabla u|^2
+2\mathbb{E}\int_0^t\int_D|u|^{p+1}+\mathbb{E}\int_0^t\int_D|\sigma(x,s)|^2.
  \eess
Denote
  \bess
v(t)=\mathbb{E}\int_Du^2, \ \ \ h(t)=\mathbb{E}\int_D\left(-2|\nabla u|^2
+2|u|^{p+1}+|\sigma(x,t)|^2\right),
   \eess
then we have
   \bess
v(t)-v(0)=\int_0^th(s)ds.
    \eess
Let
   \bess
I(t)=\int_0^tv(s)ds+A,\ \ \ A\ {\rm is \ a\ positive\ constant},
   \eess
then we have $I'(t)=v(t),\ I''(t)=h(t)$. Set
   \bess
J(t)= \mathbb{E}\int_D\left(-\frac{1}{2}|\nabla u|^2
+\frac{1}{p+1}|u|^{p+1}\right).
   \eess
It\^{o} formula implies that
   \bess
&&\frac{1}{2}\mathbb{E}\int_D|\nabla u|^2-\frac{1}{2}\mathbb{E}\int_D|\nabla u_0|^2\\
&=&-\int_0^t\mathbb{E}\int_D\Delta u(\Delta u+|u|^{p-1}u)+\frac{1}{2}\int_0^t\mathbb{E}\int_D|\nabla\sigma(x,t)|^2,
   \eess
and
   \bess
&&\frac{1}{p+1}\mathbb{E}\int_D|u|^{p+1}-\frac{1}{p+1}\mathbb{E}\int_D|u_0|^{p+1}\\
&=&\int_0^t\mathbb{E}\int_D|u|^{p-1}u(\Delta u+|u|^{p-1}u)+\frac{p}{2}\int_0^t\mathbb{E}\int_D|u|^{p-1}\sigma^2(x,t).
   \eess
Therefore, we have
   \bess
J(t)=J(0)+\int_0^t\mathbb{E}\int_D(\Delta u+|u|^{p-1}u)^2-\frac{1}{2}\int_0^t\mathbb{E}\int_D|\nabla\sigma(x,s)|^2
+\frac{p}{2}\int_0^t\mathbb{E}\int_D|u|^{p-1}\sigma^2(x,s).
  \eess
By comparing $I''(t)$ and $J(t)$, we have, for $1<\delta<\frac{p+1}{2}$,
   \bess
I''(t)=h(t)\geq4(1+\delta)J(t).
   \eess
Clearly,
   \bess
I'(t)=v(t)&=&v(0)+\int_0^th(s)ds\\
&=&v(0)+\int_0^t\mathbb{E}\int_D|\sigma(x,t)|^2+\int_0^t\mathbb{E}\int_D\left(-2|\nabla u|^2
+2|u|^{p+1}\right)dxds\\
&=&v(0)+\int_0^t\mathbb{E}\int_D|\sigma(x,t)|^2+\int_0^t\mathbb{E}\int_D\left(2u\Delta u
+2|u|^{p+1}\right)dxds.
   \eess
It follows that, for any $\varepsilon>0$,
    \bess
I'(t)^2&\leq&4(1+\varepsilon)\left[\int_0^t\mathbb{E}\int_D\left(\Delta u
+|u|^{p-1}u\right)^2dxds\right]\left[\int_0^t\mathbb{E}\int_Du^2dxds\right]\\
&&+\frac{1}{1+\varepsilon}\left[v(0)+\int_0^t\mathbb{E}\int_D|\sigma(x,t)|^2\right]^2.
   \eess
Combining the above estimates, we obtain
   \bess
&&I''(t)I(t)-(1+\alpha)I'(t)^2\\
&\geq&4(1+\delta)\left[J(0)+\int_0^t\mathbb{E}\int_D(\Delta u+|u|^{p-1}u)^2-\frac{1}{2}\int_0^t\mathbb{E}\int_D|\nabla\sigma(x,s)|^2\right.\\
&&\left.+\frac{p}{2}\int_0^t\mathbb{E}\int_D|u|^{p-1}\sigma^2(x,s)\right]
\times\left[\int_0^t\int_Du^2dxds+A\right]\\
&&-4(1+\alpha)(1+\varepsilon)\left[\int_0^t\mathbb{E}\int_D\left(\Delta u
+|u|^{p-1}u\right)^2dxds\right]\left[\int_0^t\mathbb{E}\int_Du^2dxds\right]\\
&&-\frac{(1+\alpha)}{1+\varepsilon}\left[v(0)+\int_0^t\mathbb{E}\int_D|\sigma(x,t)|^2\right]^2
  \eess
Now we choose $\varepsilon$ and $\alpha$ small enough such that
  \bess
1+\delta>(1+\alpha)(1+\varepsilon).
   \eess
By assumption,
   \bess
 J(0)-\frac{1}{2}\int_0^t\mathbb{E}\int_D|\nabla\sigma(x,s)|^2>0.
    \eess
We can choose $A$ large enough such that
   \bess
I''(t)I(t)-(1+\alpha)I'(t)^2>0,
  \eess
which implies that
   \bess
\frac{d}{dt}\left(\frac{I'(t)}{I^{1+\alpha}(t)}\right)>0.
  \eess
Then we have
   \bess
\frac{I'(t)}{I^{1+\alpha}}(t)>\frac{I'(0)}{I^{1+\alpha}(0)}\ \ \ {\rm for}\ t>0.
  \eess
It follows that $I(t)$ cannot remain finite for all $t$. This is a contradiction. The proof is complete.  $\Box$

\begin{remark}\lbl{r2.1} The advantage of concavity method is that we did not use the positivity
of solutions.
Meanwhile, the disadvantage of Theorem \ref{t2.2} is that we only deal with the additive noise.
For multiplicative noise, when we deal with the term $\mathbb{E}\int_D|\nabla u|^2$, by using
It\^{o} formula, we will have the term $-\frac{1}{2}\int_0^t\mathbb{E}\int_D|\nabla\sigma(u)|^2$, and
we cannot control this term.
   \end{remark}

\begin{remark}\lbl{r2.2} The effect of noise on the blowup problem can be described as the followings:

(i) for an additive noise, without help of the nonlinear term, the solutions will not blow up in
finite time; but if the solutions blow up in finite time without noise, the additive noise
can make the finite time blowup hard to happen. In other words,
the assumption on initial data will be stronger if we add the additive noise.

(ii) for multiplicative noise, without the help of nonlinear term, the solutions blow up
in finite time under some assumptions on initial data.
  \end{remark}

Look back at Proposition \ref{p1.2} and Theorems \ref{t2.1} and \ref{t2.2}, we find
the finite time blowup appear in the $L^p$-norm of the solutions, $p>1$. Maybe we will
ask what about the case $0<p<1$. The following result answer this equation. Consider the
following stochastic parabolic equations
        \bes\left\{\begin{array}{lll}
   du=[\Delta u+f(u)]dt+\sigma(u) dW(x,t), \ \ \qquad &t>0,\ x\in  D,\\[1.5mm]
   u(x,0)=u_0(x), \ \ \ &\qquad\ \ x\in D,\\[1.5mm]
   u(x,t)=0, \ \ \ \ &t>0,\ x\in\partial D.
    \end{array}\right.\lbl{2.5}\ees

  \begin{theo}\lbl{t2.3}
Assume $f(r)\geq0$ for $r\leq0$. Then we have:

  (i) Assume further that $f(r)\geq C_0r^p$, $q(x,y)\leq q_0$ for $x,y\in D$ and $\sigma^2(u)\leq C_1u^2$.
If the initial data satisfies
    \bess
\left(\int_Du_0(x)\phi(x)dx\right)^{p-1}>\frac{\hat\lambda}{C_0\epsilon}, \ \
\hat\lambda=\epsilon\lambda_1+\frac{\epsilon}{2}(1-\epsilon)q_0C_1^2.
  \eess
then the solution $u(x,t)$ of (\ref{2.5}) will blow up  in finite time in $L^1$-norm and $\epsilon$-order moment,
where $0<\epsilon<1$ and $p>1$, i.e., there exists a positive $T>0$ such that
   \bess
\mathbb{E}\|u(\cdot,t)\|^\epsilon_{L^1(D)}\to\infty,\ \ \ {\rm as}\ \ \ t\to T;
   \eess

(ii) Assume further that $f(r)\leq C_0r^p$, $q(x,y)\geq q_1$ for $x,y\in D$ and $\frac{1}{C_1}u^m\leq\sigma^2(u)\leq C_2u^m$.
Then, if $m>p>1$, $(m-p)(2m-1)>mp$ and the initial data are bounded,
then the solution $u(x,t)$ of (\ref{2.5}) will exist globally in the following sense:
$\mathbb{E}[|(u,\phi)|^\epsilon]<\infty$ for any $t>0$.
  \end{theo}

{\bf Proof.} (i) It follows from Proposition \ref{p1.1} that
(\ref{2.5}) has a unique positive solution. Similar to the proof of Theorem \ref{t2.1}, we will prove the theorem by contradiction.
Suppose the claim is false. Then there exists a global positive solution $u$ such
that for any $T>0$
   \bess
\sup_{0\leq t\leq T}\mathbb{E}\|u(\cdot,t)\|^\epsilon_{L^1(D)}<\infty,
   \eess
which implies that
   \bess
\sup_{0\leq t\leq T}\mathbb{E}\left(\int_ D u(x,t)\phi(x)dx\right)^\epsilon\leq\|\phi\|_{L^\infty(D)}\sup_{0\leq t\leq
T}\mathbb{E}\|u(\cdot,t)\|^\epsilon_{L^1(D)}<\infty.
    \eess
Set $\hat u=(u,\phi)$. It\^{o} formula gives that
     \bes
\hat u^\epsilon(t)&=&(u_0,\phi)^\epsilon-\epsilon\lambda_1\int_0^t\hat u^\epsilon(s)ds
+\epsilon\int_0^t\hat u(s)^{\epsilon-1}\int_ Df(u)\phi dxds\nm\\
&&+
\epsilon\int_0^t\int_ D \hat u(s)^{\epsilon-1}\sigma(u)\phi(x)d W(x,s)dx\nm\\
&&+\frac{\epsilon(\epsilon-1)}{2}\int_0^tu(s)^{\epsilon-2}\int_ D\int_ D q(x,y) \sigma(u)\phi(x)\sigma(u)\phi(y)dxdyds
   \lbl{2.6}\ees
Let $\eta(t)=\mathbb{E}\hat u^\epsilon(t)$. Similar to the proof of Theorem \ref{t2.1}, by taking an expectation over (\ref{2.6}), we obtain
  \bess
\eta(t)&=&(u_0,\phi)^\epsilon-\epsilon\lambda_1\int_0^t\eta(s)ds+
\epsilon\int_0^t\mathbb{E}\hat u(s)^{\epsilon-1}\int_ Df(u)\phi dxds\nm\\
&&+\frac{\epsilon(\epsilon-1)}{2}\int_0^t\mathbb{E}u(s)^{\epsilon-2}\int_ D\int_ D q(x,y) \sigma(u)\phi(x)\sigma(u)\phi(y)dxdyds.
    \eess
Using the assumptions $\inf_{x,y\in D}q(x,y)\leq q_0$ and $\sigma^2(u)|\leq C_1|u|^2$ and
Jensen's inequality, we have
   \bess
\eta(t)&\geq&\eta(0)-\varepsilon\lambda_1\int_0^t\eta(s)ds+C_0\epsilon\int_0^t\eta^{\frac{p+\epsilon-1}{\epsilon}}(s)ds
-\frac{\epsilon}{2}(1-\epsilon)q_0C_1^2\int_0^t\eta(s)ds,
   \eess
or, in the differential form,
   \bess\left\{\begin{array}{llll}
\ds\frac{d\eta(t)}{dt}=-\hat\lambda\eta(t)+C_0\epsilon\eta^{\frac{p+\epsilon-1}{\epsilon}}(t)\\[1.5mm]
\eta(0)=\eta_0.
    \end{array}\right. \eess
Noting that $\eta'(0)>0$. This implies that $\eta(t)>0$. An integration of the differential equation gives that
   \bess
T\leq\int_{\eta_0}^{\eta(T)}\frac{dr}{C_0\epsilon r^{\frac{p+\epsilon-1}{\epsilon}}-\hat\lambda r}
\leq\int_{\eta_0}^\infty\frac{dr}{C_0\epsilon r^{\frac{p+\epsilon-1}{\epsilon}}-\hat\lambda r}<\infty,
    \eess
which implies $\eta(t)$ must blow up at a time
$T^*\leq\int_{\eta_0}^\infty\frac{dr}{C_0\epsilon r^{\frac{p+\epsilon-1}{\epsilon}}-\hat\lambda r}$.
Hence this is a contradiction. Thus we obtain the desired result.

(ii) Define
   \bess
\tau_n=\inf\{t>0, \ \ (u,\phi)^\epsilon>n\}.
  \eess
Set $\hat u=(u,\phi)$. By using It\^{o} formula, for $t\leq\tau_n$, we have
     \bes
\hat u^\epsilon(t)&=&(u_0,\phi)^\epsilon-\epsilon\lambda_1\int_0^t\hat u^\epsilon(s)ds
+\epsilon\int_0^t\hat u(s)^{\epsilon-1}\int_ Df(u)\phi dxds\nm\\
&&+
\epsilon\int_0^t\int_ D \hat u(s)^{\epsilon-1}\sigma(u)\phi(x)d W(x,s)dx\nm\\
&&+\frac{\epsilon(\epsilon-1)}{2}\int_0^tu(s)^{\epsilon-2}\int_ D\int_ D q(x,y) \sigma(u)\phi(x)\sigma(u)\phi(y)dxdyds
   \lbl{2.7}\ees
Let $\eta(t)=\mathbb{E}\hat u^\epsilon(t)$. By taking an expectation over (\ref{2.7}), we obtain
  \bes
\eta(t)&=&(u_0,\phi)^\epsilon-\epsilon\lambda_1\int_0^t\eta(s)ds+
\epsilon\int_0^t\mathbb{E}\hat u(s)^{\epsilon-1}\int_ Df(u)\phi dxds\nm\\
&&+\frac{\epsilon(\epsilon-1)}{2}\int_0^t\mathbb{E}u(s)^{\epsilon-2}\int_ D\int_ D q(x,y) \sigma(u)\phi(x)\sigma(u)\phi(y)dxdyds
\nm\\
&\leq&\eta(0)-\epsilon\lambda_1\int_0^t\eta(s)ds+C_0
\epsilon\int_0^t\mathbb{E}\hat u(s)^{\epsilon-1}\int_ D|u|^p\phi dxds\nm\\
&&+\frac{\epsilon(\epsilon-1)}{2}\int_0^t\mathbb{E}u(s)^{\epsilon-2}\int_ D\int_ D q(x,y) \sigma(u)\phi(x)\sigma(u)\phi(y)dxdyds.
    \lbl{2.8}\ees
H$\ddot{o}$lder inequality and $\varepsilon$-Young inequality yield that
   \bess
  &&C_0\epsilon\hat u(s)^{\epsilon-1}\int_ D|u|^p\phi dx\\
  &\leq&  C_0\epsilon\hat u(s)^{\epsilon-1}\left(\int_ D|u|^m\phi dx\right)^{\frac{p}{m}}\\
  &\leq&\frac{\epsilon q_1(1-\epsilon)}{4C_1}u(s)^{\epsilon-2}\left(\int_ D|u|^m\phi dx\right)^2+C
    u(s)^{\frac{2m}{2m-p}(2p-p\epsilon-1+\epsilon)}.
      \eess
Submitting the above inequality into (\ref{2.8}), and using the assumptions on $\sigma$, we have
  \bes
\eta(t)
&\leq&\eta(0)-\epsilon\lambda_1\int_0^t\eta(s)ds+C\int_0^tu(s)^{\frac{2m}{2m-p}(2p-p\epsilon-1+\epsilon)}ds\nm\\
&&
-\int_0^t\frac{\epsilon q_1(1-\epsilon)}{2C_1}u(s)^{\epsilon-2}\left(\int_ D|u|^m\phi dx\right)^2ds\nm\\
&\leq&\eta(0)-\epsilon\lambda_1\int_0^t\eta(s)ds+C\int_0^tu(s)^{\frac{2m}{2m-p}(2p-p\epsilon-1+\epsilon)}ds\nm\\
&&
-\int_0^t\frac{\epsilon q_1(1-\epsilon)}{2C_1}u(s)^{2m+\epsilon-2}ds.
    \lbl{2.9} \ees
The assumption $(m-p)(2m-1)>mp$ gives
   \bess
\epsilon<\frac{2m}{2m-p}(2p-p\epsilon-1+\epsilon)<2m+\epsilon-2.
   \eess
Noting that for any $r<m<n$ and $u>0$, we have
  \bes
u^m=u^\beta u^{m-\beta}\leq \varepsilon u^n+C(\varepsilon)u^r, \ \ \ \beta=\frac{r(n-m)}{n-r}.
    \lbl{2.10}\ees
So we can use (\ref{2.10}) to deal with the second last term of right hand side of (\ref{2.9}).
Eventually, we get for $t\leq\tau_n$
  \bess
\eta(t)\leq\eta(0)+C\int_0^t\eta(s)ds.
    \lbl{2.11} \eess
We remark the constant $C$ does not depend on $t$. The Gronwall's lemma implies that
   \bess
\eta(t)\leq C+Ce^{Ct}, \ \ \ t\leq\tau_n.
   \eess
Letting $n\to\infty$,  the above inequality implies that $\mathbb{P}\{\tau_\infty<\infty\}=0$.
The proof is complete. $\Box$

\section{Whole space}
\setcounter{equation}{0}

In this section, we consider stochastic parabolic equations in whole space.
Our aim is to establish the global existence and non-existence under some assumptions.
We first recall the results of Foondun et al. \cite{FLN2018}, where the
authors considered the following equation
   \bes
\partial_tu_t(x)=\mathcal{L}u_t(x)+\sigma(u_t(x))\dot{F}(x,t)
 \ t>0,\ x\in\mathbb{R}^d.
   \lbl{3.1} \ees
Here $\mathcal{L}$ denotes the fractional Laplacian, the generator of an
$\alpha$-stable process and $\dot{F}$ is the random forcing term which they took to be white in time and possibly colored in space. They obtained the
following results.
   \begin{prop}\lbl{p3.1}\cite[Theorems 1.2,1,5,1.6,1.8,1.9]{FLN2018}

(i) Noise white both in time and space, i.e.,
   \bess
\mathbb{E}[\dot{F}(x,t)\dot{F}(y,s)]=\delta_0(t-s)\delta_0(x-y).
   \eess
Assume that  there exists a $\gamma>0$ such that
   \bess
\sigma(x)\geq|x|^{1+\gamma}\ \ \ {\rm for\ all\ }\ x\in\mathbb{R}^d,
  \eess
and that  there is a positive constant $\kappa$ such that $\inf_{x\in\mathbb{R}^d}:=\kappa$ .
Then
there exists a $t_0>0$ such that for all $x\in\mathbb{R}^d$, the solution
$u_t(x)$ of (\ref{3.1}) blows up in finite time, i.e.,
   \bes
\mathbb{E}|u_t(x)|^2=\infty\ \ \ {\rm whenever}\ \ t\geq t_0.
    \lbl{3.2}\ees
Furthermore, the initial condition can be weaken as the following,
    \bes
\int_{B(0,1)}u_0(x)dx:=K_{u_0}>0,
   \lbl{3.3}\ees
where $B(0,1)$ is the ball centred in the point $0$ and radius $1$.
The solution
$u_t(x)$ of (\ref{3.1}) also blows up in finite time whenever $K_{u_0}\geq K$, where $K$ is some positive constant.

(ii) Noise white in time and correlated in space, i.e.,
   \bess
\mathbb{E}[\dot{F}(x,t)\dot{F}(y,s)]=\delta_0(t-s)f(x,y).
   \eess
Assume that for fixed $R>0$, there exists some positive number $K_f$ such
that
     \bes
\inf_{x,y\in B(0,R)}(x,y\in B(0,R))f(x,y)\geq K_f.
  \lbl{3.4}\ees
Then, for fixed $t_0>0$ there exists a positive unmber $\kappa_0$ such
that for all $\kappa\geq\kappa_0$ and $x\in\mathbb{R}^d$ we have (\ref{3.2}) holds.

In particularly, suppose that the correlation function $f$ is
given by
    \bess
f(x,y)=\frac{1}{|x-y|^\beta} \ \ {\rm with}\ \ \beta<\alpha\wedge d.
  \eess
Then for $\kappa>0$ there exists a $t_0>0$ such that (\ref{3.2}) holds.

Furthermore, under the assumptions (3.3) and (\ref{3.4}), there
exists a $t_0>0$ such that for all $x\in\mathbb{R}^d$ (\ref{3.2}) holds.
  \end{prop}
In the above proposition, Foondun et al. \cite{FLN2018} only considered
the finite time blowup phenomenon driven by noise. Our aim in this paper
is to find the effect of noise, including additive noise and multiplicative noise. And we are also very interested in the type (3) as introduction
said.

We first consider the global existence of the following stochastic parabolic equations
   \bes \left\{\begin{array}{llll}
du_t=(\Delta u+f(u,x,t))dt+\sigma(u,x,t)dB_t,\ \ t>0,\ &x\in\mathbb{R}^d,\\
u(x,0)=u_0(x)\gneqq0, \ \  &&x\in\mathbb{R}^d,
   \end{array}\right.\lbl{3.5}\ees
where $B_t$ is one-dimensional Brownian motion.
A mild solution to (\ref{3.5}) in sense of Walsh \cite{walsh1986} is any
$u$ which is adapted to the filtration generated by the
white noise and satisfies the following evolution equation
   \bess
u(x,t)&=&\int_{\mathbb{R}^d}K(t,x-y)u_0(y)dy
+\int_0^t\int_{\mathbb{R}^d}K(t-s,x-y)f(u,y,s)dyds\\
&&
+\int_0^t\int_{\mathbb{R}^d}K(t-s,x-y)\sigma(u,y,s)dydB_s,
   \eess
where $K(t,x)$ denotes the heat kernel of Laplacian operator, i.e.,
   \bess
K(t,x)=\frac{1}{(2\pi t)^{d/2}}\exp\left(-\frac{|x|^2}{2t} \right)
  \eess
satisfies
   \bess
\left(\frac{\partial}{\partial t}-\Delta\right)K(t,x)=0\ \ \ {\rm for}\ \
(x,t)\neq(0,0).
   \eess
We get the following results.
   \begin{theo}\lbl{t3.1}
Suppose that there exist positive constants $C_0,\ 0<p<1$ such that
   \bess
|h(u,x,t)|\leq C_0|u|^p,\ \ \
h=f \ {\rm or}\ g.
  \eess
Then the solutions of (\ref{3.5}) with bounded continuous
initial data $u_0$ exist globally in any
$r$-order moment, $r\geq1$.
  \end{theo}

{\bf Proof.} By taking the second moment and using the Walsh isometry,
we get for any $T>0$
   \bess
\mathbb{E}|u(x,t)|^2&=&\left(\int_{\mathbb{R}^d}K(t,x-y)u_0(y)dy
+\int_0^t\int_{\mathbb{R}^d}K(t-s,x-y)f(u,y,s)dyds\right.\\
&&\left.
+\int_0^t\int_{\mathbb{R}^d}K(t-s,x-y)\sigma(u,y,s)dydB_s)\right)^2\\
&\leq&4\int_{\mathbb{R}^d}K(t,x-y)u^2_0(y)dy+4C_0^2
\int_0^t\int_{\mathbb{R}^d}K(t-s,x-y) [\mathbb{E}|u(y,s)|^2]^pdyds\\
&&
+4C_0^2\int_0^t\mathbb{E}\left(\int_{\mathbb{R}^d}K(t-s,x-y) |u(y,s)|^p]\right)^2\\
&\leq&4 \sup_{x\in\mathbb{R}^d}|u_0(x)|^2+8C_0^2\sup_{t\in[0,T],x\in\mathbb{R}^d}[\mathbb{E}|u(y,s)|^2]^p
\int_0^T\int_{\mathbb{R}^d}K(t,x) dtdx
   \eess
Then taking supremum for $t,x$ over $\in[0,T]\times\mathbb{R}^d$ (the right hand is independent of $t$ and $x$) , we get
    \bess
\sup_{t\in[0,T],x\in\mathbb{R}^d} \mathbb{E}|u(x,t)|^2\leq   4 \sup_{x\in\mathbb{R}^d}|u_0(x)|^2+8C_0^2T\sup_{t\in[0,T],x\in\mathbb{R}^d}[\mathbb{E}|u(y,s)|^2]^p.
   \eess
Notice that $0<p<1$, we have for any $T>0$
   \bess
\sup_{t\in[0,T],x\in\mathbb{R}^d} \mathbb{E}|u(x,t)|^2\leq C(T)<\infty,
  \eess
which implies that $\mathbb{P}\{|u(x,t)|=\infty\}=0$. The proof is complete. $\Box$

We remark that the heat kernel $K$ belongs to $L^1(\mathbb{R}^d)$ but not $L^2(\mathbb{R}^d)$.
Hence this result does not hold for the noise white in both time and space. Meanwhile, if we assume
the covariance function $q(x,y)$ is uniformly bounded, then the above result also hold for
the noise white in time and  correlated in space.

Next, we establish the result similar to the case of type (3). In order to do that,
we will consider the following Cauchy problem
   \bes \left\{\begin{array}{llll}
du_t=\Delta udt+\sigma(u,x,t)dW(x,t),\ \ t>0,\ &x\in\mathbb{R}^d,\\
u(x,0)=u_0(x)\gneqq0, \ \  &&x\in\mathbb{R}^d,
   \end{array}\right.\lbl{3.6}\ees
where $W(t,x)$ is white noise both in time and space.
In the rest of paper, we always assume that the initial data is nonnegative continuous function.
A mild solution to (\ref{3.6}) in sense of Walsh \cite{walsh1986} is any
$u$ which is adapted to the filtration generated by the
white noise and satisfies the following evolution equation
   \bess
u(x,t)=\int_{\mathbb{R}^d}K(t,x-y)u_0(y)dy
+\int_0^t\int_{\mathbb{R}^d}K(t-s,x-y)\sigma(u,y,s)W(dy,ds),
   \eess
where $K(t,x)$ denotes the heat kernel of Laplacian operator.
We get the following results.
   \begin{theo}\lbl{t3.2}
Suppose $d=1$ and $\sigma^2(u,x,t)\geq C_0u^{2m}$, $C_0>0$, then for $1<m\leq\frac{3}{2}$, the solutions of (\ref{3.5}) blows up in finite
time for any nontrivial nonnegative initial data $u_0$. That is to say, there exists a positive constant $T$ such
that for all $x\in\mathbb{R}$
   \bess
\mathbb{E}u^2(x,t)=\infty\ \ {\rm for}\ t\geq T.
  \eess
  \end{theo}

{\bf Proof.} We assume that the solution remains finite for all finite $t$ almost
surely and want to derive a contradiction.
By taking the second moment and using the Walsh isometry,
we get
   \bess
\mathbb{E}|u(x,t)|^2&=&\left(\int_{\mathbb{R}^d}K(t,x-y)u_0(y)dy\right)^2
+\int_0^t\int_{\mathbb{R}^d}K^2(t-s,x-y)\mathbb{E}\sigma^2(u,y,s)dyds\\
&=:&I^2_1(x,t)+I_2(x,t).
   \eess
We may assume without loss of generality that $u_0(x)\geq C_1>0$ for $|x|<1$ by the
assumption. A direct computation shows that
   \bes
I_1(x,t)&\geq&\frac{C_1}{(2\pi t)^{d/2}}\int_{B_1(0)}\exp\left(-\frac{|x|^2+|y|^2}{2t}\right)dy\nm\\
&\geq& \frac{C_1}{(2\pi t)^{d/2}}\exp\left(-\frac{|x|^2}{2t}\right)\int_{|y|\leq\frac{1}{\sqrt{t}}}
\exp\left(-\frac{|y|^2}{2}\right)dy\nm\\
&\geq& \frac{C}{(2\pi t)^{d/2}}\exp\left(-\frac{|x|^2}{2t}\right)
  \lbl{3.7} \ees
for $t>1$ and $C>0$.

It is easy to see that
   \bess
I_2(x,t)&\geq& C_0\int_0^t\int_{\mathbb{R}^d}K^2(t-s,x-y)\mathbb{E}|u(y,s)|^{2m}dyds\\
&\geq& C_0\int_0^t\int_{\mathbb{R}^d}K^2(t-s,x-y)[\mathbb{E}|u(y,s)|^{2}]^mdyds.
   \eess
Denote $v(x,t)=\mathbb{E}|u(x,t)|^{2}$. Let
   \bess
G(t)=\int_{\mathbb{R}^d}K(t,x) v(x,t)dx.
   \eess
Then for $t>1$,
   \bes
G(t)&=&\int_{\mathbb{R}^d}I^2_1(x,t)K(t,x)dx+\int_{\mathbb{R}^d}I_2(x,t)K(t,x)dx\nm\\
&\geq&\frac{C_2}{t^d}+\int_0^t\int_{\mathbb{R}^d}\int_{\mathbb{R}^d}K(t,x)K^2(t-s,x-y)v^m(y,s)dydxds.
   \lbl{3.8}\ees
It is clear that
   \bess
&&\int_{\mathbb{R}^d}K(t,x)K^2(t-s,x-y)dx\\
&=&\frac{1}{(2\pi t)^{d/2}[2\pi (t-s)]^{d}}\int_{\mathbb{R}^d}
\exp\left(-\frac{|x|^2}{2t}-\frac{|x-y|^2}{t-s}\right)dx\\
&=&K(s,y)\frac{(2\pi s)^{d/2}}{(2\pi t)^{d/2}[2\pi (t-s)]^{d}}\int_{\mathbb{R}^d}
\exp\left(\frac{|y|^2}{2s}-\frac{|x|^2}{2t}-\frac{|x-y|^2}{t-s}\right)dx.
   \eess
Since
   \bess
&&\frac{|y|^2}{2s}-\frac{|x|^2}{2t}-\frac{|x-y|^2}{t-s}\\
&\geq&\frac{|y|^2}{2s}-\frac{|x-y|^2+|y|^2+2|x-y||y|}{2t}-\frac{|x-y|^2}{t-s}\\
&=&\frac{1}{2t}\left(-2|x-y||y|+\frac{t-s}{s}|y|^2\right)-\frac{|x-y|^2}{2t}-\frac{|x-y|^2}{t-s}\\
&\geq&-\frac{s|x-y|^2}{2t(t-s)}-\frac{|x-y|^2}{2t}-\frac{|x-y|^2}{t-s}\\
&\geq& -\frac{2|x-y|^2}{t-s}\ \ {\rm for}\ 0<s<t,
  \eess
we get for $0<s<t$
   \bess
\int_{\mathbb{R}^d}
\exp\left(\frac{|y|^2}{2s}-\frac{|x|^2}{2t}-\frac{|x-y|^2}{t-s}\right)dx
\geq \int_{\mathbb{R}^d}
\exp\left(-\frac{2|x-y|^2}{t-s}\right)dx
=C_3(t-s)^{d/2}.
  \eess
Substituting the above estimate into (\ref{3.8}) and applying Jensen's inequality, we obtain
   \bess
G(t)&\geq&\frac{C_2}{t^d}+C_4\int_0^t\frac{s^{d/2}}{t^d}\int_{\mathbb{R}^d}K(s,y)v^m(y,s)dydxds\\
&\geq&\frac{C_2}{t^d}+C_4\int_0^t\frac{s^{d/2}}{t^d}G^m(s)ds
   \eess
We can rewrite the above inequality as
   \bes
t^d G(t)&\geq& C_2+C_4\int_0^ts^{d/2}G^m(s)ds\lbl{3.9}\\
&=:&g(t).\nm
   \ees
Then for $t>1$, we have
   \bess
&&g(t)\geq C_2,\\
&&g'(t)\geq C_4t^{d/2}G^m(t)\geq C_4t^{d/2}\left(\frac{1}{t^d}g(t)\right)^m=C_4t^{\frac{d}{2}-dm}g^m(t),
  \eess
which implies
   \bess
\frac{C_2^{1-m}}{m-1}\geq \frac{1}{m-1}g^{1-m}(t)\geq C_4\int_t^Ts^{\frac{d}{2}-dm}dx\ \ {\rm for}\
T>t\geq1.
   \eess
If $m\leq\frac{d+2}{2d}$, that is, $\frac{d}{2}-dm+1\geq0$, the right-hand side of the above
inequality is unbounded as $T\to\infty$, which gives a contradiction. Noting that we must
let $m>1$ because we used the Jensen's inequality, thus we get $1<m\leq\frac{3}{2}$ and $d=1$.
And thus we complete
the proof. $\Box$

If the noise is just one-dimensional Brownian motion, the result
will be different. For this, we consider the following stochastic
  \bes \left\{\begin{array}{llll}
du_t=\Delta udt+\sigma(u,x,t)dB_t,\ \ t>0,\ &x\in\mathbb{R}^d,\\
u(x,0)=u_0(x)\gneqq0, \ \  &&x\in\mathbb{R}^d,
   \end{array}\right.\lbl{3.10}\ees
where $B_t$ is one-dimensional Brownian motion.
A mild solution to (\ref{3.10}) in sense of Walsh \cite{walsh1986} is any
$u$ which is adapted to the filtration generated by the
white noise and satisfies the following evolution equation
   \bess
u(x,t)=\int_{\mathbb{R}^d}K(t,x-y)u_0(y)dy
+\int_0^t\int_{\mathbb{R}^d}K(t-s,x-y)\sigma(u,y,s)dydB_s,
   \eess
where $K(t,x)$ denotes the heat kernel of Laplacian operator.

   \begin{theo}\lbl{t3.3}
Suppose $d=1$ and $\sigma^2(u,x,t)\geq C_0u^{2}$, $C_0>0$, then the solutions of (\ref{3.10}) blows up in finite
time for any nontrivial nonnegative initial data $u_0$.
  \end{theo}

{\bf Proof.} Similar to the proof of Theorem \ref{t3.2}, we assume that the solution remains finite for all finite $t$ almost
surely.
By taking the second moment and using the Walsh isometry,
we get
   \bess
\mathbb{E}|u(x,t)|^2&=&\left(\int_{\mathbb{R}^d}K(t,x-y)u_0(y)dy\right)^2
+\int_0^t\left(\int_{\mathbb{R}^d}K(t-s,x-y)\mathbb{E}\sigma(u,y,s)dy\right)^2ds\\
&=:&u_1(x,t)+u_2(x,t).
   \eess
We may assume without loss of generality that $u_0(x)\geq C_1>0$ for $|x|<1$ by the
assumption.
The estimate (\ref{3.7}) also holds, i.e.,
   \bess
u_1(x,t)\geq \frac{C}{(2\pi t)^{d/2}}\exp\left(-\frac{|x|^2}{2t}\right)
    \eess
for $t>1$ and $C>0$.

It is easy to see that, for $m\geq2$,
   \bess
u_2(x,t)&\geq& C_0\int_0^t\left(\int_{\mathbb{R}^d}K(t-s,x-y)\mathbb{E}|u(y,x)|^{m}dy\right)^2ds\\
&\geq& C_0\int_0^t\left(\int_{\mathbb{R}^d}K(t-s,x-y)[\mathbb{E}|u(y,x)|^{2}]^{m/2}dy\right)^2ds\\
&\geq& C_0\int_0^t\left(\int_{\mathbb{R}^d}K(t-s,x-y)\mathbb{E}|u(y,x)|^{2}dy\right)^mds.
   \eess
Denote $v(x,t)=\mathbb{E}|u(x,t)|^{2}$. Let
   \bess
G(t)=\int_{\mathbb{R}^d}K(t,x) v(x,t)dx.
   \eess
Then for $t>1$,
   \bes
G(t)&=&\int_{\mathbb{R}^d}u^2_1(x,t)K(t,x)dx+\int_{\mathbb{R}^d}u_2(x,t)K(t,x)dx\nm\\
&\geq&\frac{C_2}{t^d}+\int_0^t\left(\int_{\mathbb{R}^d}\int_{\mathbb{R}^d}K(t,x)K(t-s,x-y)v(y,s)dydx\right)^mds.
   \lbl{3.11}\ees
It is clear that (see \cite[Page 42]{Hubook2018})
   \bess
\int_{\mathbb{R}^d}K(t,x)K(t-s,x-y)dx\geq C_3K(s,y)\left(\frac{s}{t}\right)^{d/2}.
   \eess
Substituting the above estimate into (\ref{3.11}) and applying Jensen's inequality, we obtain
   \bess
G(t)\geq\frac{C_2}{t^d}+C_3\int_0^t\left(\frac{s^{d/2}}{t^{d/2}}\right)^mG^m(s)ds
   \eess
We can rewrite the above inequality as
   \bes
t^{md/2} G(t)&\geq& C_2t^{(m-2)d/2}+C_3\int_0^ts^{dm/2}G^m(s)ds\lbl{3.9}\\
&=:&g(t).\nm
   \ees
Then for $t>1$, we have
   \bess
&&g(t)\geq C_2t^{(m-2)d/2},\\
&&g'(t)\geq C_3t^{dm/2}G^m(t)\geq C_3t^{d/2}\left(\frac{1}{t^{dm/2}}g(t)\right)^m=C_3t^{(1-m)md/2}g^m(t),
  \eess
which implies
   \bess
\frac{C_2^{1-m}}{m-1}t^{-d(m-1)(m-2)/2}\geq \frac{1}{m-1}g^{1-m}(t)\geq C_4\int_t^Ts^{(1-m)md/2}dx\ \ {\rm for}\
T>t\geq1.
   \eess
If $(m-1)md/2\leq1$, we will get a contradiction by letting $T\to\infty$. If
$\frac{d(m-1)(m-2)}{2}>-1+ \frac{(m-1)md}{2}$, then we will get a contradiction by letting
$T\to\infty$ and then taking $t\gg1$. Noting that when $m=2,\ d=1$, we have
$(m-1)md/2=1$ and $\frac{d(m-1)(m-2)}{2}>-1+ \frac{(m-1)md}{2}$ is equivalent to $m<1+\frac{1}{d}$.
Since $m\geq2$, we get a contradiction for the case that $m=2,\ d=1$. The proof is complete. $\Box$

\begin{remark} \lbl{r3.1}
Comparing Theorem \ref{t3.2} with Proposition \ref{p3.1},
the assumptions of Proposition \ref{p3.1} on initial data need the lower bound,
but in Theorem \ref{t3.2} we did not.

Theorems \ref{t3.2} and \ref{t3.3} show that the time-space white noise and
Brownian motion are different. But the method used here
is not suitable to fractional Laplacian operator.
Sugitani \cite{Sug1975} established the Fujita index for
Cauchy problem of fractional Laplacian operator.
The main difficult is that we can not get the exact estimate of
$\int_{\mathbb{R}^d}p^2(t,x)p(t-s,x-y)dx$, where $p(t,x)$ is the
heat kernel of fractional Laplacian operator.
  \end{remark}

\section{Discussion}
\setcounter{equation}{0}
An interesting issue of stochastic partial differential equations is to find
the difference when we add the noise, i.e., the impact of
noise. For stochastic partial differential equations, we want to know whether the solutions keep positive.
In this section, we first consider the positivity of the solutions of stochastic parabolic
equations in the whole space, and then consider the impact of noise.

In the followings, we will select a test function $\beta_\varepsilon(r)$.
Define
    \bess
&&\beta_\varepsilon(r)=\int_r^\infty\rho_\varepsilon(s)ds,\ \ \
\rho_\varepsilon(r)=\int_{r+\varepsilon}^\infty J_\varepsilon(s)ds,\ \ \ r\in\mathbb{R},\nm\\
&&J_\varepsilon(|x|)=\varepsilon^{-n}J\left(\frac{|x|}{\varepsilon}\right),\ \ \
J(x)=\left\{\begin{array}{llll}\medskip
C\exp\left(\frac{1}{|x|^2-1}\right),\ \ \ &|x|<1,\\
0, \ \ \ &|x|\geq1.
  \end{array}\right. \eess
Then by direct verification, we have the following result.

\begin{lem}\lbl{l4.1} The above constructed functions $\rho_\varepsilon,\beta_\varepsilon$ are in $C^\infty(\mathbb{R})$ and
have the following properties:
$\rho_\varepsilon$ is a non-increasing function and
   \bess
\beta_\varepsilon'(r)=-\rho_\varepsilon(r)=\left\{\begin{array}{llll}\medskip
0,\ \ \ &r\geq0,\\
-1,\ \ \ &r\leq-2\varepsilon.
 \end{array}\right.\eess
 Additionally,  $\beta_\varepsilon$ is convex and
    \bess
\beta_\varepsilon(r)=\left\{\begin{array}{llll}\medskip
0,\ \  \ \ \ &r\geq0,\\
-2\varepsilon-r+\varepsilon \hat C, \ \  &r\leq-2\varepsilon,
   \end{array}\right.\eess
where $\hat C=\int^0_{-2}\int_{t+1}^1J(s)dsdt<2$. Furthermore,
   \bess
0\leq\beta_\varepsilon''(r)=J_\varepsilon(r+\varepsilon)\leq \varepsilon^{-d}C, \ \ -2\varepsilon\leq r\leq0,
   \eess
which implies that
    \bess
-2^dC\leq r^d\beta_\varepsilon''(r)\leq0 \ \ \ \ & {\rm for}\ -2\varepsilon\leq r\leq0,\ {\rm and}\ d\ {\rm is}\ {\rm odd};\\
0\leq r^d\beta_\varepsilon''(r)\leq2^dC\ \ \ \ &{\rm for}\ -2\varepsilon\leq r\leq0,\ {\rm and}\ d\ {\rm is}\ {\rm even}.
   \eess
   \end{lem}

Now, we consider the following stochastic parabolic  equations
    \bes\left\{\begin{array}{lll}
   du=(\Delta u+f(u,x,t))dt+g(u,x,t)dW(x,t), \ \ \qquad t>0,\ &x\in \mathbb{R},\\[1.5mm]
   u(x,0)=u_0(x), \ \ \ \ \ &x\in \mathbb{R},
    \end{array}\right.\lbl{4.1}\ees
where $W(x,t)$ is time-space white noise.

\begin{theo}\lbl{t4.1} Assume that (i) the function $f(r,x,t)$ is continuous on
  $\mathbb{R}\times\mathbb{R}\times[0,T]$;  (ii) $f(r,x,t)\geq0$ for $r\leq0$, $x\in\mathbb{R}$
and $t\in[0,T]$; and (iii)  $ g^2(u,x,t)\leq ku^{2m}$,
where $k>0$, $2m>1$ and $(-1)^{2m-1}\in\mathbb{R}$. Then the solution of
initial-boundary value problem {\rm(\ref{4.1})} with nonnegative initial datum remains positive:
$u(x,t)\geq0$, a.s. for almost every $x\in \mathbb{R}$ and for all $  t\in[0,T]$.
\end{theo}

{\bf Proof.}  Define
   \bess
\Phi_\varepsilon(u_t)=(1,\beta_\varepsilon(u_t))=\int_{\mathbb{R}} \beta_\varepsilon(u(x,t))dx.
   \eess
By It\^{o}'s formula, we have
    \bess
\Phi_\varepsilon(u_t)&=&\Phi_\varepsilon(u_0)+\int_0^t\int_{\mathbb{R}} \beta_\varepsilon'(u(x,s))\Delta u(x,s)dxds\\
&&+\int_0^t\int_ {\mathbb{R}} \beta_\varepsilon'(u(x,s))f(u(x,s),x,s)dxds\\
&&+\int_0^t\int_ {\mathbb{R}} \beta_\varepsilon'(u(x,s))g(u(x,s),x,s)dW(x,s)dx\\
&&+\frac{1}{2}\int_0^t\int_{\mathbb{R}^d}\beta_\varepsilon''(u(x,s))g^2(u(x,s),x,t)dxds\\
&=&\Phi_\varepsilon(u_0)+\int_0^t\int_{\mathbb{R}}\beta_\varepsilon''(u(x,s))\left(\frac{1}{2}g^2(u(x,s),x,s)-|\nabla u|^2\right)dxds\\
&&+\int_0^t\int_ {\mathbb{R}} \beta_\varepsilon'(u(x,s))f(u(x,s),x,s)dxds\\
&&+\int_0^t\int_ {\mathbb{R}} \beta_\varepsilon'(u(x,s))g(u(x,s),x,s)dW(x,s)dx.
   \eess
Taking expectation over the above equality and using Lemma \ref{l4.1}, we get
    \bess
\mathbb{E}\Phi_\varepsilon(u_t)
&=&\mathbb{E}\Phi_\varepsilon(u_0)+\mathbb{E}\int_0^t\int_{\mathbb{R}} \beta_\varepsilon''(u(x,s))\\
&&\times\left(\frac{1}{2}g^2(u(x,s),x,s)-|\nabla u|^2\right)dxds\\
&&+\mathbb{E}\int_0^t\int_{\mathbb{R}} \beta_\varepsilon'(u(x,s))f(u(x,s),x,s)dxds\\
&\leq&\mathbb{E}\Phi_\varepsilon(u_0)+\frac{k}{2}\mathbb{E}\int_0^t\int_{\mathbb{R}} \beta_\varepsilon''(u(x,s))
u(x,s)^{2m}dxds\\
&&+\mathbb{E}\int_0^t\int_{\mathbb{R}}\beta_\varepsilon'(u(x,s))f(u(x,s),x,s)dxds.
  \eess
Here and after, we denote $\|\cdot\|_{L^1}$ by $\|\cdot\|_1$. Let $\eta(u)=u^-$
denote the negative part of $u$ for $u\in\mathbb{R}$.
Then we have $\lim\limits_{\varepsilon\rightarrow0}\mathbb{E}\Phi_\varepsilon(u_t)=\mathbb{E}\|\eta(u_t)\|_1$.
It follows from Lemma \ref{l4.1} that
    \bess
0\geq u^{2m}\beta''_\varepsilon(u)\geq\left\{\begin{array}{llll}\medskip
0,\ \ \ \ & u\geq0\ {\rm or}\ u\leq-2\varepsilon,\\
-2Cu^{2m-1},\ \ \ \ &-2\varepsilon\leq u\leq0,\ {\rm and}\ u^{2m-1}\geq0,
   \end{array}\right.
   \eess
or
    \bess
0\leq u^{2m}\beta''_\varepsilon(u)\leq\left\{\begin{array}{llll}\medskip
0,\ \ \ \ & u\geq0\ {\rm or}\ u\leq-2\varepsilon,\\
-2Cu^{2m-1},\ \ \ \ &-2\varepsilon\leq u\leq0,\ {\rm and}\ u^{2m-1}\leq0
   \end{array}\right.
   \eess
which implies that
$\lim\limits_{\varepsilon\rightarrow0}u^{2m}\beta_\varepsilon''(u)=0$ provided that $2m>1$.
By taking the limits termwise as $\varepsilon\rightarrow0$ and using Lemma \ref{l4.1}, we get
    \bess
\mathbb{E}\|\eta(u_t)\|_1&\leq&\mathbb{E}\|\eta(u_0)\|_1
-\mathbb{E}\int_0^t\int_{\mathbb{R}} \eta'(u(x,s))f(u(x,s),x,s)dxds\nm\\
&\leq&0,
   \eess
which implies that $u^-=0$ a.s. for a.e. $x\in D$, $\forall t\in[0,T]$.
This completes the proof.
$\Box$

If $W(x,t)$ is replaced by $B_t$ in (\ref{4.1}), then Theorem \ref{t4.1} holds for any dimension.
The reason why we only consider one dimension in Theorem \ref{t4.1} is that the It\^{o} formula
only holds for one-dimensional time-space white noise.

In order to find the impact of noise, we first recall a well-known result of
deterministic parabolic equations. Consider the Cauchy problem
      \bes \left\{\begin{array}{llll}
\frac{\partial}{\partial t}u_t=\Delta u+u^p,\ \ t>0,\ &x\in\mathbb{R}^d,\\
u(x,0)=u_0(x)\gneqq0, \ \  &x\in\mathbb{R}^d.
   \end{array}\right.\lbl{4.2}\ees
\begin{prop}\lbl{p4.1} (i) If $p>1+\frac{2}{d}$, then the solution of
(\ref{4.2}) is global in time, provided the initial datum satisfies,
for some small $\varepsilon>0$,
  \bess
u_0(x)\leq\varepsilon K(1,x), \ \ \ x\in\mathbb{R}^d.
  \eess
(ii) If $1<p\leq1+\frac{2}{d}$, then all nontrivial solutions of
(\ref{4.2}) blow up in finite time.
  \end{prop}
Next we consider the stochastic parabolic equation
      \bes \left\{\begin{array}{llll}
du_t=[\Delta u+|u|^{p-1}u]dt+\sigma(u)dW(x,t),\ \ t>0,\ &x\in\mathbb{R}^d,\\
u(x,0)=u_0(x)\gneqq0, \ \  &x\in\mathbb{R}^d.
   \end{array}\right.\lbl{4.3}\ees
It is well known that the mild solution of (\ref{4.3})  can be written as
   \bess
u(x,t)&=&\int_{\mathbb{R}^d}K(t,x-y)u_0(y)dy+\int_0^t\int_{\mathbb{R}^d}K(t-s,x-y)|u|^{p-1}udyds\\
&&
+\int_0^t\int_{\mathbb{R}^d}K(t-s,x-y)\sigma(u,y,s)W(dy,ds).
   \eess
 \begin{theo}\lbl{t4.2} Assume all the assumptions of Theorem
\ref{t4.1} hold. Then $1<p\leq1+\frac{2}{d}$, then the expectation of all nontrivial solutions of
(\ref{4.3}) blow up in finite time. That is to say, there exists a positive constant $t_0>0$
such that $\mathbb{E}u(x,t)=\infty,\  \ t\geq t_0$ for all $x\in\mathbb{R}^d$. When $m>1$, the mean square of
solutions to (\ref{4.3}) will blow up in finite time under the
condition that the initial data is suitable large.
  \end{theo}

{\bf Proof.} It follows from Theorem \ref{t4.1} that the solutions of
(\ref{4.3}) keep positive. Following the representation of mild solution, we have
  \bess
\mathbb{E}u(x,t)=\int_{\mathbb{R}^d}K(t,x-y)\mathbb{E}u_0(y)dy+\int_0^t\int_{\mathbb{R}^d}K(t-s,x-y)\mathbb{E}|u|^pdyds,
   \eess
which implies that
   \bess
\mathbb{E}u(x,t)\geq\int_{\mathbb{R}^d}K(t,x-y)\mathbb{E}u_0(y)dy+\int_0^t\int_{\mathbb{R}^d}K(t-s,x-y)[\mathbb{E}u]^pdyds,
   \eess
Denoting $v(x,t)=\mathbb{E}u(x,t)$, we have that $v(x,t)$ is a super-solution of (\ref{4.2}).
By the results of Proposition \ref{p4.1} and comparison principle, we obtain that there exists a positive constant $t_0>0$
such that $\mathbb{E}u(x,t)=\infty,\  \ t\geq t_0$ for all $x\in\mathbb{R}^d$. Meanwhile,
noting that
   \bess
\mathbb{E}u(x,t)\leq\left(\mathbb{E}u^p(x,t)\right)^{\frac{1}{p}},\ \ p>1,
  \eess
we have that $\mathbb{E}u^p(x,t)$, $p>1$, will blow up in finite time.

When $m>1$, we have
   \bess
\mathbb{E}|u(x,t)|^2&\geq&\left(\int_{\mathbb{R}^d}K(t,x-y)u_0(y)dy\right)^2
+\int_0^t\int_{\mathbb{R}^d}K^2(t-s,x-y)\mathbb{E}\sigma^2(u,y,s)dyds\\
&=:&w(x,t).
   \eess
Foondun \cite{FLN2018} proved the mean square of function $w(x,t)$ will blow up in finite
time under the condition that the initial data is suitable large. So the solution $u$ will
also blow up in finite time.
The proof is complete. $\Box$

\medskip

\noindent {\bf Acknowledgment} The first author was supported in part
by NSFC of China grants 11771123. The authors thanks Prof. Feng-yu Wang for discussing
this manuscript.


\begin{thebibliography}{99}
\bibitem{BY2014} J. Bao and C. Yuan, {\em Blow-up for stochastic reactin-diffusion equations with jumps},
J Theor. Probab. {\bf29} (2016) 617-631.

\bibitem{Cb2007} P-L. Chow, {\em Stochastic partial differential equations},
Chapman Hall/CRC Applied Mathematics and Nonlinear Science Series. Chapman Hall/CRC,
Boca Raton, FL, 2007. x+281 pp. ISBN: 978-1-58488-443-9.

\bibitem{C2009} P-L. Chow, {\em  Unbounded positive solutions of nonlinear parabolic It\^{o}
equations}, Communications on Stochastic Analysis {\bf 3} (2009) 211-222.

\bibitem{C2011} P-L. Chow, {\em  Explosive solutions of stochastic reaction-diffusion
equations in mean $L^p$-norm}, J. Differential Equations {\bf 250} (2011) 2567-2580.


\bibitem{CL2012} Pao-Liu Chow and K. Liu, {\em  Positivity and explosion in mean $L^p$-norm
of stochastic functional parabolic equations of retarded type}, Stochastic Processes and their
Applications {\bf 122} (2012) 1709-1729.

\bibitem{DMH2015}A. Debussche, S. de Moor and M. Hofmanov$\acute{a}$, {\em A regularity result for quasilinear stochastic
partial differential equations of parabolic type}, SIAM J. Math. Anal.  {\bf47}  (2015) 1590-1614.

\bibitem{DW2014} J. Duan and W. Wang, {\em Effective Dynamics of Stochastic
Partial Differential Equations}, Elsevier, 2014.

\bibitem{DL2010} M. Dozzi and J. A. L\'{o}pez-Mimbela, {\em
Finite-time blowup and existence of global positive solutions of a semi-linear spde},
Stochastic Process. Appl., {\bf120} (2010) 767-776.

\bibitem{FF2013} E. Fedrizzi and F. Flandoli, {\em Noise prevents singularities in linear transport equations},
Journal of Functional Analysis {\bf264} (2013) 1329-1354.

\bibitem{FLN2018} M. Foondun, W. Liu and E. Nane {\em Some non-existence
results for a class of stochastic partial differential equations}, J. Differential Equations in press.

\bibitem{F1966} H. Fujita, {\em On the blowing up of solutions of the Cauchy problem for
$u_t-\Delta u=u^{1+\alpha}$}, J. Fac. Sci. Univ. Tokyo Sect.
IA Math. {\bf 13} (1966) 109-124.

\bibitem{F1970} H. Fujita, {\em On some nonexistence and nonuniqueness theorems for nonlinear parabolic equations},
 Proc. Symp. Pure Math. {\bf XVIII} (1970) 105-113.


\bibitem{Hubook2018} B. Hu, {\em Blow-up Theories for Semilinear Parabolic Equations},
Lecture Notes in Mathematics ISSN print edition: 0075-8434, Springer Heidelberg Dordrecht London New York, 2018.

\bibitem{kH1973} K. Hayakawa, {\em On nonexistence of global solutions of some semilinear parabolic equations}, Proc. Japan Acad. Ser.
A Math. {\bf49} (1973) 503-505.


\bibitem{LPJ2016} K. Li, J. Peng and J. Jia, {\em Explosive solutions
of parabolic stochastic partial differential equations with L$\acute{e}$vy
noise}, arXiv:1306.01676.


\bibitem{LR2010} W. Liu and M. R\"{o}ckner, {\em SPDE in Hilbert space with locally monotone
coefficients}, J. of Functional Analysis {\bf 259} (2010) 2902-2922.


\bibitem{L2013} W. Liu, {\em Well-posedness of stochastic partial differential equations with
Lyapunov condition}, J. Differential Equations {\bf 254} (2013) 725-755.

\bibitem{LD2015} G. Lv and J.   Duan, {\em Impacts of Noise on a Class of Partial Differential Equations},
J. Differential Equations, {\bf258} (2015) 2196-2220.

\bibitem{LWW2016}G. Lv, L. Wang and X. Wang, {\em Positive and unbounded solution of stochastic delayed evolution equations}, Stoch. Anal. Appl. {\bf34} (2016) 927-939.

\bibitem{MZ1999} R. Manthey and T. Zausinger, {\em Stochastic evolution equations in $L^{2\nu}_\rho$},
Stochastics and Stochastic Report {\bf 66} (1999) 37-65.

\bibitem{M1991} C. Mueller, {\em Long time existence for the heat equation with a noise term},
Probab. Theory Related Fields {\bf 90} (1991) 505-517.

\bibitem{MuS1993} C. Mueller and R. Sowers, {\em Blowup for the heat equation with a noise term},
Probab. Theory Related Fields {\bf 93} (1993) 287-320.

\bibitem{NX2012} M. Niu and B. Xin, {\em Impacts of Gaussian noises on the blow-up times of nonlinear stochastic
partial differential equations}, Nonlinear Analysis: Real World Applications {\bf13} (2012) 1346-1352.

\bibitem{eP1979} E. Pardoux, {\em Stochastic partial differential equations and filtering of diffusion processes}, Stochastic, {\bf3} (1979) 127-167.

\bibitem{PZ1992} G. Da Prato and J. Zabczyk, {\em Stochastic equations in infinite dimensions},
Encyclopedia of Mathematics and its applications, Cambridge University Press (1992).

\bibitem{PZ1992} G. Da Prato and J. Zabczyk, {\em Nonexplosion, boundedness and ergodicity
for stochastic semilinear equations}, J. Differential Equations {\bf 98} (1992) 181-195.


\bibitem{SGKM} A. Samarskii, V. Galaktionov, S. Kurdyumov and S. Mikhailov, {\em
Blow-up in quasilinear parabolic equations}, Walter de Gruyter, Berlin, New York, 1995.


\bibitem{Shiga} T. Shiga {\em Some properties of solutions for one-dimensional SPDE's associated
with space-time white noise}, Gaussian random fields (Nagoya, 1990), 354-363.

\bibitem{Sug1975} S. Sugitani, {\em On nonexistence of global
solutions for some nonlinear integral equations}, Osaka J. Math.,
{\bf12} (1975) 35-51.


\bibitem{T2009} T. Taniguchi, {\em The existence and uniqueness of energy solutions to local
non-Lipschitz stochastic evolution equations}, J. Math. Anal. Appl. {\bf 360} (2009) 245-253.


\bibitem{walsh1986} John B. Walsh, {\em
An introduction to Stochastic Partial Differential Equations},
volume 1180 of Lecture Notes in Math., pages 265-439, Springer
Berlin, 1986.


\end{thebibliography}
 \end{document}